\documentclass[11pt]{article}
\usepackage{amsfonts}
\usepackage{amsmath, amsthm}
\usepackage{amssymb}
\usepackage{eucal}

\begin{document}

\baselineskip=16pt
\title{A counterexample for the conjecture of finite simple
groups}

\author{{Wujie Shi}\\
{\small Faculty of Mathematics and Finance, Chongqing University of
Arts and Sciences}\\
{\small  Chongqing 402160, P. R. China}\\
{and}\\
{\small School of Mathematics, Suzhou University}\\
{\small Suzhou, 215006, P. R. China}\\
{\small\  E-mail: shiwujie@outlook.com}}
\date{}
\maketitle
\begin{abstract} In this note we provide some counterexamples for the
conjectures of finite simple groups, one of the conjectures said
"all finite simple groups $G$ can be determined using their orders
$|G|$ and the number of elements of order $p$, where $p$ the largest prime divisor of $|G|$".\\


\end{abstract}

For a finite group, the order of group and the element orders are
two of the most important basic concepts. Let $G$ be a finite group
and $\pi_e(G)$ be the set of element orders of
$G$. In 1987, the author posed the following conjecture[1]:\\

{\bf Conjecture.} Let $G$ be a group and $S$ a finite simple group.
Then $G\cong S$ if and only if (a) $\pi_e(G) =  \pi_e(S)$, and (b) $|G| = |S|$.\\

That is, for all finite simple groups we may characterize them using
only their orders and the sets of their element orders (briefly,
"two orders").\\

From 1987 to 2003, the authors of [1-7] proved that this conjecture
is correct for all finite simple groups except $B_n$, $C_n$ and
$D_n$ ($n$ even). In 2009, the authors of [8] proved that this
conjecture is correct for $B_n$, $C_n$ and $D_n$ ($n$ even). Thus,
this conjecture is proved and become a theorem, that is, {\bf all
finite simple groups can determined by their
"two orders"}.\\

{\bf Question 1.}  Weaken the condition of "two orders",
characterize all finite simple groups.\\

Assume that $|G| = p_1^{\alpha_1} p_2^{\alpha_2}...p_t^{\alpha_t}$,
$p_1 < p_2 < ... < p_t$ and $\pi_e(G) = \{1, ... , m_3, m_2, m_1\}$,
$1 < ... < m_3 < m_2 < m_1$. Whether we may use the part of "order
of group"(i.e. the order of Hall-subgroup of $G$) or the maximum
part of $\pi_e(G)$(i.e. $m_1$, $\{m_2, m_1\}$, or $\{m_3, m_2,
m_1\}$) to characterize all finite simple groups?

{\bf Remark.} A remarkable result of [9] shows the three largest
element orders $\{m_3, m_2, m_1\}$ determine the characteristic of
simple groups of Lie type of odd characteristic.\\

{\bf Result $1.^{[10, Theorem 14]}$} Let $G$ be a group and $p>3$ a
prime. Suppose that $|G| = |PSL(2, p)| = (p - 1)p(p + 1)/2$ and
$m_1(G) = m_1(PSL(2, p)) = p, p \neq 7$. Then $G \cong PSL(2, p)$.
For $p =7$, if $|G|= 168$ and $\{m_2, m_1\} = \{4, 7\}$, then $G
\cong
PSL(2, 7)$.\\

{\bf Question 2.} Compare to $\pi_e(G)$, consider other
(conjugated)invariant quantitative sets of $G$ and discuss their influence on $G$.\\

Let $G$ be a group and $nse(G) = \{m_k |k \in \pi_e(G)\}$, where
$m_k$ denote the number of elements of order $k$ in $G$.\\

{\bf Result $2.^{[11, Theorem]}$} A group $G$ (which is not assumed
finite) is isomorphic to $A_i$ if and only if $nse(G) = nse(A_i), i = 4,5,6$.\\

{\bf Result $3.^{[12, Main Theorem]}$} Let $G$ be a group and $S$ a
simple $K_4$-group (i.e. $|\pi(S)|=4$). Then $G \cong S$ if and only
if the following hold: (1) $|G| = |S|$,(2) $nse(G) = nse(S)$.\\

{\bf Result $4.^{[13]}$} If $G$ is a group such that $nse(G) =
nse(PSL(2, q))$, where $q \in \{7, 8, 11, 13\}$, then $G \cong PSL(2, q)$.\\

In [14], the author took a new perspective and show that some
families of simple groups are basically determined just by the
number of elements of order $p$, where $p$ is the largest prime
divisor of the order of the group. Then the author posed the
following two conjectures:\\

{\bf Conjecture C.} Let $S$ be a finite simple group and $p$ the
largest prime divisor of $|S|$. If $G$ is a finite group with the
same number of elements of order $p$ as $S$ and $|G| = |S|$, then $G
\cong S$.\\

{\bf Conjecture E.} Let $S$ be a finite simple group that is not
isomorphic to $L_2(q)$, where $q$ is a Mersenne prime and let $p$ be
the largest prime divisor of $|S|$. If a finite group $G$ is
generated by elements of order $p$ and $G$ has the same number of
elements of order $p$ as $S$, then $G/Z(G) \cong S$.\\

In this note we provide some counterexamples for the above
conjectures. We construct the counterexample of Conjecture C
firstly.

{\bf Counterexample 1.} Let $S = A_8 \cong L_4(2)$. Then $|S| = 2^6.
3^2.5.7$ and $p = 7$ is the largest prime. Since $\pi_e(A_8) = \{1,
2, 3, ... ,7, 15\}$,
we have $|C_S(<a>)| = 7$, where $o(a) = 7$.\\

From $N_S(<a>)/C_S(<a>) \cong$ a subgroup of $Aut(<a>)$ and
$|S:N_S(<a>)| \equiv 1(mod\,\,7)$, we have $S$ contains $2^7.3^2.5$
elements of order 7.\\

Let $G = L_3(4)$. From $|G| = 2^6. 3^2.5.7$, $p = 7$ and
 $\pi_e(L_3(4)) = \{1, 2, 3, 4, 5 ,7\}$, we have $G$ contains $2^7.3^2.5$
elements of order 7 similarly. $L_3(4)$ is not isomorphic to $A_8$.\\

{\bf Counterexample 2.} Let $S = O_7(3)$ and $G = S_6(3)$. We have
$|S| = |G| = 2^9.3^9.5.7.13$ and $p = 13$ is the largest prime.
Since $\pi_e(O_7(3)) = \{1, 2, 3, ... ,10, 12, 13, 14, 15, 18, 20\}$
and $\pi_e(S_6(3)) = \{1, 2, 3, ... ,10, 12, 13, 14, 15, 18, 20, 24,
30, 36\}$, we may give same conclusion as Counterexample 1.\\

{\bf Counterexample 3.} Let $S = L_2(7)$ and $G = [[Z_2 \times Z_2
\times Z_2]Z_7]Z_3$, a 2-Frobenius group. We have $|S| = |G| =
2^3.3.7$ and $p = 7$ is the largest prime. Since $\pi_e(L_2(7)) =
\{1, 2, 3, 4, 7\}$ we have $|C_S(<a>)| = 7$, where $o(a) = 7$. We
have $S$ contains $2^4.3$ elements of order 7. Similarly $G$
contains $2^4.3$ elements of order 7, it is also a counterexample for Conjecture C.\\

For Conjecture E, $A_8$ and $L_3(4)$ is also a counterexample. Put
$S = A_8$ and $G = L_3(4)$, then 7 is the largest prime divisor of
$|A_8|= |L_3(4)|$. We get the number of elements of order 7 is
$2^7.3^2.5$ in $S$ by calculation. The order of subgroup generated
by all elements of order 7 greater than $2^7.3^2.5$, and great than
the orders of all maximal subgroups of $A_8$(see [15]). Thus $A_8$
is generated by elements of order $7$. By the same argument,
$L_3(4)$ has the same number of elements of order $7$ as $A_8$, but
$L_3(4)$ is not isomorphic to $A_8$, we get the counterexample.
Moreover, $O_7(3)$ and $S_6(3)$ is also a counterexample for Conjecture E.\\

However, if we avoid these same order and non-isomorphic simple
groups, is the conjecture E hold? This is a question that we can be
further studied.\\

Related the number of elements of same order Prof. J.G. Thompson
posed the following open problem$^{[16]}$.\\

{\bf Definition.} Groups $G_1$ and $G_2$ are called same order type
groups, if $\pi_e(G_1) = \pi_e(G_2)$, and for all $k \in
\pi_e(G_1)$, the number of elements of order $k$ in $G_1$ and $G_2$
are equal. In particular, $nse(G_1) = nse(G_2)$.\\

In 1987,  J.G. Thompson gave the following nonsolvable example of
same order type groups which are not isomorphic$^{[16]}$.\\

$G_1 = 2^4 : A_7$ and $G_2 = L_3(4) : 2_2$ both are maximal
subgroups of $M_{23}$.\\

In the above example, $|G_1| = |G_2| = 2^7.3^2.5.7$, $\pi_e(G_1) =
\pi_e(G_2) = \{1,2,3,4,5,6,7,8,14\}$ and $nse(G_1) = nse(G_2)=\{1,
435, 2240, 6300, 8064, 6720, 5040, 5760\}$. $p = 7$ is the largest
prime and the number of elements of order $p$ is 5760.\\

"This example shows that there are finite groups $G_1, G_2$ of the
same order type which do not have the same set of composition
factors"$^{[16]}$\\

{\bf Open Problem.$^{[16]}$} Suppose $G_1, G_2$ are finite groups of
same order type. Suppose also that $G_1$ is solvable. Is it true
that $G_2$ is also necessarily solvable?\\

\end{document}